\documentclass[draft]{amsart}
\usepackage{amsmath,amsthm,amssymb}
\newtheorem{lemma}{Lemma}
\newtheorem{corollary}{Corollary}
\newtheorem{proposition}{Proposition}

\newcommand{\C}{\mathbb{C}\mkern1mu}
\newcommand{\R}{\mathbb{R}\mkern1mu}

\newcommand{\bmat}{\begin{pmatrix}}
\newcommand{\emat}{\end{pmatrix}}
\newcommand{\bsmat}{\bigl(\begin{smallmatrix}}
\newcommand{\esmat}{\end{smallmatrix}\bigr)}
\newcommand{\Bsmat}{\Bigl(\begin{smallmatrix}}
\newcommand{\Esmat}{\end{smallmatrix}\Bigr)}
\newcommand{\bbsmat}{\biggl(\begin{smallmatrix}}
\newcommand{\eesmat}{\end{smallmatrix}\biggr)}

\newcommand{\gen}[1]{\langle #1\rangle}

\DeclareMathOperator{\Ad}{Ad}
\DeclareMathOperator{\ad}{ad}

\newcommand{\mf}{\mathfrak}
\begin{document}
\title{Biquotient actions on unipotent Lie groups}
\author{Annett P\"uttmann}
\date{}
\thanks{Research supported by SFB/TR12 ``Symmetrien und Universalit\"at in mesoskopischen 
Systemen'' of the Deutsche Forschungsgemeinschaft.}
\address{Fakult\"at f\"ur Mathematik, Ruhr-Universit\"at Bochum, 44780 Bochum}
\email{annett.puettmann@rub.de}
\begin{abstract}
We consider pairs $(V,H)$ of  subgroups of a connected unipotent  
complex Lie group $G$ for which the induced $V\times H$-action on $G$ 
by multiplication from the left and from the right is free.
We prove that this action is proper if the Lie algebra $\mf{g}$ of $G$  is $3$-step nilpotent.
If $\mf{g}$ is $2$-step nilpotent, then there is a global slice of the action that is isomorphic
to $\C^n$.
Furthermore,  a global slice isomorphic to $\C^n$ exists if $\dim V = 1 = \dim H$ or 
$\dim V = 1$ and $\mf{g}$ is $3$-step nilpotent.
We give an explicit example of a $3$-step nilpotent Lie group and a pair of 
$2$-dimensional subgroups such that the
induced action is proper but the corresponding geometric quotient is not affine.
\end{abstract}
\maketitle
\section{Introduction and generic situation}
By $G$ we always denote a connected, simply connected complex unipotent Lie group, 
i.e. an algebraic subgroup of a group of upper triangular matrices.
Any pair $(V,H)$ of complex algebraic subgroups $V,H\subset G$ defines
a natural $V\times H$ action on  $G$ by $(v,h).g=vgh^{-1}$.
Equivalently, we can discuss the $V$-action given by $v.gH=vgH$
on the homogeneous space $G/H$, 
which is isomorphic to $\C^{\dim G-\dim H}$ as a complex algebraic variety. 
We consider only free $V \times H$-actions on $G$.

Denote by $\mf{g}$, $\mf{v}$, and $\mf{h}$ the Lie algebras of $G$, $V$, and $H$, respectively.
We will frequently use the descending central series
$\mf{g} = \mf{g}^{(0)} \supset \ldots \supset \mf{g}^{(j+1)}:=[\mf{g},\mf{g}^{(j)}] \supset
	\ldots \supset \{ 0\}$. 
A nilpotent Lie algebra $\mf{g}$ is said to be $l$-step nilpotent if $l$ is the smallest integer
such that $\mf{g}^{(l)} =\{ 0\}$.

Since the only compact algebraic subgroup of $V$ is the trivial group and the exponential map
$\exp: \mf{g} \to G$ is an isomorphism,
the isotropy groups $V_x$ are compact for all $x \in G/H$ iff
the isotropy groups $V_x$ are trivial for all $x\in G/H$ iff
$\Ad(g)(\mf{v}) \cap \mf{h} =\{ 0\}$ for all $g\in G$.

Lipsman conjectured \cite{Lip} that the $V\times H$-action on $G$ is proper if it is free.
In fact, Nasrin proved \cite{Nasrin} Lipsman's conjecture to be true for 
$2$-step nilpotent Lie groups $G$.
We improve this result showing that Lipsman's conjecture is true if $\mf{g}^{(3)} = \{ 0\}$, 
and that there is a global slice isomorphic to $\C^{\dim G -\dim H -\dim V}$, 
if $\mf{g}^{(2)} = \{ 0\}$.

There is a counterexample to Lipsman's conjecture presented by Yoshino \cite{YT1}
with $\mf{g}^{(4)} = \{ 0\}$. 
Essentially the same example serves to show that there is a free non-proper affine 
$\C^2$-action on $\C^5$, which is just the smallest member of a series, $n\geq 5$, 
of free non-proper affine $\C^2$-actions on $\C^n$ \cite{AP}.
We construct a $3$-step nilpotent Lie algebra $\mf{g}$ and $2$-dimensional
subalgebras $\mf{h}$ and $\mf{v}$, such that the induced $V$-action on $G/H$ is
Winkelmann's \cite{Wink}
the free affine proper $\C^2$-action on $\C^6$ without global slice.

We can not answer the question of  whether Lipsman's conjecture is true 
if one of the subgroups is one-dimensional. 
But a global slice isomorphic to $\C^n$ exists if $\dim V = 1 = \dim H$ or 
$\dim V = 1$ and $\mf{g}$ is $3$-step nilpotent.

Let us start by looking at the generic case which turns out to be easy:
We call a basis ${\mathcal B}=((X_1,\ldots,X_n))$ of $\mf{g}$ a Levi-Malcev basis, 
if it is compatible with the descending central series, i.e.,
${\mathcal B}\cap \mf{g}^{(j)}$ is a basis of $\mf{g}^{(j)}$ for all $j$.
Given subalgebras $\mf{v},\mf{h}\subset\mf{g}$
a vector space decomposition $\mf{g} = \mf{v}\oplus \mf{s} \oplus \mf{h}$ is called
Levi-Malcev decomposition if  there are bases of the subspaces 
$\mf{v}$, $\mf{h}$, and $\mf{s}$ such that their union is a Levi-Malcev basis of $\mf{g}$.
Note that there is a Levi-Malcev decomposition if and only if
$\pi_{j}(\mf{v}\cap\mf{g}^{(j)})\cap\pi_j(\mf{h}\cap\mf{g}^{(j)}) = \{ 0\} $ for all $j$, 
where $\pi_j:\mf{g}^{(j)}\to \mf{g}^{(j)}/\mf{g}^{(j+1)}$ denotes the canonical homomorphism
from $\mf{g}^{(j)}$ to the commutative Lie algebra $ \mf{g}^{(j)}/\mf{g}^{(j+1)}$.
If two subspaces $\mf{h}, \mf{v} \subset \mf{g}$ are generically chosen, then there exists
a Levi-Malcev decomposition $\mf{g} = \mf{v} \oplus \mf{s} \oplus \mf{h}$.

\begin{lemma}
If there is a Levi-Malcev decomposition $\mf{g} = \mf{v} \oplus \mf{s} \oplus \mf{h}$, 
then $\Ad(g)(\mf{v})\cap\mf{h} = \{ 0\}$ for all $g\in G$.
\end{lemma}
\begin{proof}
If $X \in \mf{g}^{(j)}$, then $\Ad(g)(X)-X \in \mf{g}^{(j+1)}$ for all $g\in G$.
\end{proof}
Since the exponential map $\exp: \mf{g} \to G$ is an isomorphism, the action
can be pulled-back onto the Lie algebra.
We use the notation $\exp(X*Y) = (\exp X) (\exp Y)$ for $X,Y \in \mf{g}$.
\begin{proposition}
If there is a Levi-Malcev decomposition  $\mf{g} = \mf{v}\oplus \mf{s} \oplus \mf{h}$,
then $S:=\exp\mf{s} \subset G$ is a global slice of the $V\times H$-action on $G$.
\end{proposition}
\begin{proof}
If $X\not \in \mf{g}^{(j+1)}$, then $ \pi_{j} (X*Y) = \pi_{j} (X+Y) $ for all $Y \in \mf{g}$.
\end{proof}
\begin{corollary}
There is a global slice of the $H$-action on $G$.
In particular, $G/H \cong \C^{\dim G-\dim H}$.
\end{corollary}
\begin{proof}
Arbitrary  vector space decompositions 
$\pi_j(\mf{g}^{(j)}) = \mf{s}_j \oplus \pi_j(\mf{h}\cap\mf{g}^{(j)})$
imply a Levi-Malcev decomposition
$\mf{g} = \mf{s}\oplus\mf{h}$.
Then $S:=\exp\mf{s} \subset G$ is a global slice of the $H$ action on $G$.
\end{proof}
\section{Reductions using normal subgroups}
\label{sec:reductions}
If $N\lhd G$ is a normal subgroup, then $G/N$ is again a connected simply connected
unipotent Lie group. Furthermore, there is the induced action of the subgroups
$V/(V\cap N)$ and $H/(H\cap N)$ on $G/N$.  
In particular, if $Z(G)$ is the center of $G$, then $V\cap Z(G)$ and $H\cap Z(G)$ are central
subgroups of $G$.
Therefore, the $V\times H$-action on $G$ is equivalent to the 
$V/(V\cap Z(G)) \times H/(H\cap Z(G))$-action on $G/((Z(G)\cap V) \times (Z(G)\cap H))$.
Consequently, we can assume $V\cap Z(G) = H\cap Z(G) = \{ e\}$ and 
$\mf{v}\cap\mf{g}^{(l-1)} = \mf{h}\cap\mf{g}^{(l-1)} = \{ 0 \}$, which might be a weaker condition.
\begin{lemma}
Let $N\lhd G$ be a normal subgroup of $G$.
If the induced $V\times H$-action on $G/N$ is free, then any local slice $S\subset G/N$
gives a local slice $SN\subset G$.
\end{lemma}
\begin{proof}
If there are $v\in V$, $h\in H$, $s,s' \in S$, and $n,n' \in N$ such that
$vsnh^{-1} = s'n'$, then there  exists $n'' \in N$ such that $vsh n'' = s'$, 
since $N$ is normal.
\end{proof}

Let $G_1 \lhd G$ be a normal subgroup that contains $V$ and $H$.
Any global slice $S$ of the $G_1$-action on $G$ by left-multiplication defines an isomorphism
$\Psi:S \times G_1 \to G$, $(s,g_1)\mapsto sg_1$.
Now, the $V\times H$-action on $G$ can be regarded as a  family of 
$V\times H$-actions on $G_1$ parameterized by $s\in S$:
\[ (v,h).(s,g_1) = \Psi^{-1}(vsg_1h^{-1}) = \Psi^{-1}(s(s^{-1}vs)g_1h^{-1}) = (s,s^{-1}vsg_1h^{-1}) . \]
\begin{lemma}
\label{lem:family}
Let $G_1 \lhd G$ be a normal subgroup that contains $H$ and $V$.
If there is an element $Y_0 \in \mf{g}$ such that $\ad(Y_0)(\mf{v}) \subset \mf{v}$ and
$\mf{g} = \gen{Y_0}_\C \oplus \mf{g}_1$ is a Levi-Malcev decomposition,
then the $V\times H$-action on $G$ is proper iff the $V \times H$-action on $G_1$ is
proper. In that case, $V\backslash G /H \cong (V\backslash G_1 /H) \times \C$.
\end{lemma}
\begin{proof}
The condition $\ad(Y_0)(\mf{v}) \subset \mf{v}$ implies
$\Ad(\exp(y_0Y_0))(\mf{v}) = e^{y_0\ad(Y_0)}(\mf{v}) = \mf{v}$ for all $y_0 \in \C$, 
i.e., we have a trivial family of $V\times H$-actions on $G_1$.
\end{proof}

\begin{corollary}
\label{cor:reduceG1}
Let $\mf{n} := \pi_0^{-1}(\pi_0(\mf{h}) \cap \pi_0(\mf{v}))$, $\mf{v}_0 := \mf{n} \cap \mf{v}$,
$\mf{h}_0 := \mf{n} \cap \mf{h}$, 
$V_0 :=\exp \mf{v}_0$, and $H_0 := \exp \mf{h}_0$.
Let $\mf{g}_1\subset \mf{g}$ be a maximal subspace that satisfies $\mf{n} \subset \mf{g}_1$,  
$\mf{g}_1 \cap \mf{h} = \mf{h}_0$, and  $\mf{g}_1 \cap \mf{v} = \mf{v}_0$.
The $V\times H$-action on $G$ is proper iff the $V_0 \times H_0$-action on $G$ is proper.
In that case, $V\backslash G /H \cong V_0 \backslash (\exp\mf{g}_1) /H_0$.
\end{corollary}
\begin{proof}
The subspace $\mf{g}_1$ is an ideal, since it contains $\mf{g}^{(1)}$.
Note that $V$ and $H$ can be interchanged in Lemma \ref{lem:family}.
We therefore assume $\pi_0(\mf{h}) \subset \pi_0(\mf{v})$.
Choose Lie algebra elements $X_1, \ldots, X_m \in \mf{v}$ 
such that $\dim \mf{s} = m$ and $\mf{g} = \mf{s} \oplus \mf{g}_1$ 
is a Levi-Malcev decomposition for subspace $\mf{s} :=  \gen{X_1, \ldots, X_m}_\C$.
Now, $\{ \exp(x_1 X_1) \ldots \exp (x_m X_m) : x_j \in \C\}=: S$ is a global slice for the
$G_1$-action on $G$ be left multiplication.
\end{proof}

\begin{lemma}
\label{lem:X0maximal}
Let $N \lhd G$ be a normal subgroup. 
If $S_N\subset G/N$ is a global slice of the $V/(V\cap N) \times H/(H\cap N)$-action on 
$G/N$ and $S \subset G$ is a global slice of the 
$(V\cap N)\times (H\cap N)$-action on $G$,
then $S \cap S_NN$ is global slice of the $V\times H$-action on $G$.
\end{lemma}
\begin{proof}
The map $V \times (S \cap S_NN) \times H \to G$, $(v,s,g)\mapsto vsg$, is an isomorphism,
because $(V\cap N)S_NN (H\cap N) = S_N N$ and $S/N = G/N$.
\end{proof}
\section{Special cases} 
If $G$ is commutative, then any vector space decomposition 
$\mf{g} = \mf{v} \oplus \mf{s} \oplus \mf{h}$ is a Levi-Malcev decomposition
and defines a global slice $S:=\exp\mf{s}$ of the $V\times H$-action.
\subsection{$2$-step nilpotent Lie algebras}
Applying the results of section \ref{sec:reductions} we can assume
that there are $X_1, \ldots, X_m  \in \mf{g}\setminus\mf{g}^{(1)}$ and 
$Z_1,\ldots,Z_m \in \mf{g}^{(1)}$ such that
$(( X_1,\ldots, X_m))$ is a basis of $\mf{h}$ and 
$(( X_1+Z_1,\ldots, X_m+Z_m))$ is a basis of $\mf{v}$ if $\mf{g}$ is $2$-step nilpotent.
In particular, $\mf{v}$, $\mf{h}$, and $\mf{h} \oplus \mf{g}^{(1)}$ are commutative.
\begin{proposition}
If $\mf{g}$ is $2$-step nilpotent, then there exists
a global slice $S\subset G$ of the $V\times H$-action on $G$ that is
algebraically isomorphic to $\C^{\dim G -\dim H - \dim V}$.
\end{proposition}
\begin{proof}
Let $\mf{g} = \mf{s}\oplus \mf{h}$ be any Levi-Malcev decomposition.
Then $\exp\mf{s}$ is a global slice of the $H$-action on $G$.
The induced $V$-action on $\mf{s}\cong G/H$ is given by
\begin{align*} 
\exp(\sum_{j=1}^m t_j(X_j+Z_j)).Y = &
	\exp(\sum_{j=1}^m t_j(X_j+Z_j)) \exp(Y) \exp(\sum_{j=1}^m -t_jX_j) \\
= & \Ad(\exp(\sum_{j=1}^m t_jX_j))(Y)+\sum_{j=1}^m t_jZ_j \\
= & \exp(\ad(\sum_{j=1}^m t_jX_j))(Y)+\sum_{j=1}^m t_jZ_j 
\end{align*}
which is an affine action of the commutative group $\C^m$ on $\mf{s} \cong \C^{\dim G -m}$
of degree one, i.e., the expression $\exp(\sum_{j=1}^m t_j(X_j+Z_j)).Y$ is linear in the
variables $t_j$.
By \cite{AP} this action has a global slice that is algebraically isomorphic to $\C^{\dim\mf{s}-m}$.
\end{proof}
\subsection{$3$-step nilpotent Lie algebras}
If $\mf{g}$ is $3$-step nilpotent, we can assume that 
$\mf{v}\cap \mf{g}^{(2)} = \mf{h} \cap \mf{g}^{(2)}  = \{ 0\}$ and $\pi_0(\mf{v})=\pi_0(\mf{h})$.
Choose $X_1,\ldots, X_m \in \mf{h}\setminus \mf{g}^{(1)}$ such that
$((\pi_0(X_1),\ldots, \pi_0(X_m) ))$ is a basis of $\pi_0(\mf{h}) = \pi_0(\mf{v})$.
There are $Z_1,\ldots, Z_m \in \mf{g}^{(1)}$ such that $X_j+Z_j \in \mf{v}$ for $1\leq j \leq m$.
Let us consider the subalgebras $\mf{h}_0 \subset \mf{h}$ and $\mf{v}_0 \subset \mf{v}$
generated by the elements $X_j \in \mf{h}$, $j=1,\ldots,m$, 
and $X_j+Z_j \in \mf{v}$, $j=1,\ldots,m$, respectively.
Note that $[\mf{h}_0,\mf{h}\cap \mf{g}^{(1)}] = [\mf{v}_0,\mf{v}\cap \mf{g}^{(1)}] = \{ 0\}$, 
since these commutators are in $\mf{g}^{(2)}$ and $V$ and $H$ are assumed to have 
trivial intersection with the center of $G$.
Hence, there is a Lie algebra decomposition $\mf{v} = \mf{v}_0\oplus \mf{v}_1$ with
$\mf{v}_1 \subset \mf{g}^{(1)}$.

By construction there is an isomorphism $\phi: \mf{v}_0\to \mf{h}_0$, $X_j+Z_j \mapsto X_j$.
The set $\mf{n} := \{ (Y_0+Y_1,\phi(Y_0)+Y) : Y_j \in \mf{v}_j, Y \in \mf{h}\cap \mf{g}^{(1)}\} 
\subset \mf{v} \oplus \mf{h}$ is an ideal isomorphic to $\mf{v} \oplus \mf{h}\cap \mf{g}^{(1)}$.
The resulting action of the normal subgroup $\exp \mf{n} \subset V \times H$ 
is affine of degree two, 
since $X*Y = X+Y+\frac{1}{2}[X,Y]$ if $X \in \mf{g}^{(1)} = \mf{g}^{(l-2)}$, i.e.,
$\delta^3(x_j) = 0$ for all $j$ and for any locally nilpotent derivation $\delta \in \mf{n}$
for any choice of affine coordinates $x_j$ on $\mf{g}$.
\begin{lemma}
A free affine action of a connected simply connected complex unipotent Lie group $N$ on $\C^m$ of degree two is proper.
\end{lemma}
\begin{proof}
Since the action is affine, there are coordinates $x_1,\ldots,x_m$ of $\C^m$ such that
the action map $N\times \C^m \to \C^m$, $(n,x)\mapsto n.x$, defines a representation
of the Lie algebra $\mf{n}$ on $\C[x_1,\ldots,x_m]$ by locally nilpotent triangular derivations
such that $\delta(x_k) \in \gen{1,x_1,\ldots,x_{k-1}}_\C$ 
for all $\delta\in\mf{n}$ and all $k=1,\ldots,m$.
Given a pair $(x_0,\delta_0) \in \C^m\times\mf{n}$ we can assume 
$x_j(x_0) = 0$ for all $j=1,\ldots,m$ and
$\delta_0(x_k) = c_0$ and $|c_0| \ne 0$ for some $k\in \{ 1,\ldots,m\}$, since the action is free.
Consequently, there is a neighborhood $U_0$ of $(x_0,\delta_0)$ in $\C^m\times \mf{n}$ such
that $|\delta(x_k)(x)| > \frac{3}{4}c_0 $ for all $(x,\delta) \in U_0$.

We have to show that the map $\Phi : N\times \C^m \to \C^m \times \C^m$, 
$(n,x)\mapsto (x,n.x)$, is proper, i.e., for any compact set $K\subset \C^m\times \C^m$ 
the preimage $\Phi^{-1}(K)$ is compact.
Since the action is of degree two,
\[ x_j(\exp(\delta).x) = (\exp(-\delta).x_j)(x) = x_j(x)-\delta(x_j)(x) +\tfrac{1}{2}\delta^2(x_j)(x) \]
for all $\delta\in\mf{n}$, $x\in\C^m$, $j = 1,\ldots, m$.
Note that if $|(x,n.x)| \leq R$ and $n = \exp \delta$, then 
$h_j(\delta,x) := |-\delta(x_j)(x) +\frac{1}{2}\delta^2(x_j)(x)| \leq 2R$ for all $j=1,\ldots,m$.
We will prove that the set 
$\{ (x,\delta) \in \C^m\times \mf{n} :|x|<R, h_j(x,\delta) \leq 2R\,\forall j \}$ is bounded. 
For any point $(x_0,\delta_0) \in \C^m\times \mf{n}$ satisfying $|x_0|\leq R$ and $|\delta_0| =1$ 
we will find a neighborhood $U_0$ and a number $S_{(x_0,\delta_0)}$ such that
$h_j(x,\alpha\delta)\leq 2R$ for $(x,\delta) \in U_0$ and for all $j$  implies 
$|\alpha|\leq S_{(x_0,\delta_0)}$.  

For $\alpha \in \R^{\geq 0}$ and $(x,\delta) \in U_0$,
\begin{multline*} 
|-\alpha\delta(x_k)(x)+\tfrac{1}{2}\alpha^2\delta^2(x_k)(x)|  = 
	\alpha |-\delta(x_k)(x)+\tfrac{1}{2}\alpha\delta^2(x_k)(x) | \\
\geq \alpha (|\delta(x_k)(x)| - \tfrac{1}{2}\alpha |\delta^2(x_k)(x)|) 
\geq \alpha (\tfrac{3}{4}c_0 - \tfrac{1}{2}\alpha |\delta^2(x_k)(x)|).
\end{multline*}
Since $0 = \delta^3(x_k) = \sum_{j<k} a_j(\delta) \delta^2 x_j$ for all $\delta \in\mf{n}$
and $h_j(x,\alpha\delta)\leq 2R$, 
\begin{multline*}
\alpha\delta^2(x_k) = \alpha \sum_{j<k} \delta(c_0(\delta)+a_j(\delta)x_j) 
	= \sum_{j<k} a_j(\delta) \alpha\delta x_j \\
= \sum_{j<k} a_j(\delta)(\alpha\delta x_j -\tfrac{1}{2}\alpha^2\delta^2 x_j) 
	\leq 2R  \sum_{j<k} a_j(\delta).
\end{multline*}
The neighborhood $U_0$ can be shrunk further to obtain 
$\sum_{j<k} a_j(\delta) \leq \frac{c_0}{4R}$ for all $(x,\delta)\in U_0$, because
$\sum_{j<k} a_j(\delta_0) =0$. Hence, $\frac{1}{2}\alpha\delta^2(x_k)(x) \leq \frac{1}{4}c_0$.
This gives the upper bound $S_{(x_0,\delta_0)} = \frac{4R}{c_0} \geq |\alpha|$.
\end{proof}
\begin{corollary}
\label{cor:treestep}
If $\mf{g}$ is $3$-step nilpotent, then the $V\times H$-action on $G$ is proper.
\end{corollary}
\begin{proof}
By the previous Lemma the $N$-action on $G$ is proper.
This means that there are local slices $S_\alpha \subset G$ of the $N$-action such that
the sets $U_\alpha := N. S_\alpha$ cover $G$. Let $G^{(1)} := \exp \mf{g}^{(1)}$.
The $N$-action on $G/G^{(1)}$ is trivial and $(V\times H)/N \cong H/(H\cap  G^{(1)})$.
Let $\pi_0(\mf{g}) = \mf{s}_0 \oplus \pi_0(\mf{h})$ be any vector space decomposition.
Then $\tilde{U}_\alpha := \exp(\pi_0^{-1}(\mf{s}_0)) \cap U_\alpha$ are $N$-invariant subsets 
and $\tilde{S}_\alpha := \tilde{U}_\alpha /N$ are local slices of the $V\times H$-action
on $G$.
\end{proof}

To conclude this section we give an example of a $3$-step nilpotent Lie algebra $\mf{g}$
and two subalgebras $\mf{v}, \mf{h} \subset \mf{g}$ such that the
corresponding free action of $V\times H$ on $G$ does not  have 
a global (holomorphic or algebraic) slice.
It was shown in \cite{Wink} that the two commuting derivations
$\delta' = y_1\frac{\partial}{\partial y_3}+y_2\frac{\partial}{\partial y_4}
	+(1+y_1)\frac{\partial}{\partial z_2}+z_2\frac{\partial}{\partial z_1}$
and
$\delta = y_2\frac{\partial}{\partial y_3}+ \frac{\partial}{\partial z_1}$
define a proper free affine $\C^2$-action on $\C^6$ that has a geometric quotient that is neither
affine nor Stein.
\begin{lemma}
\label{lem:WinkelmannLiealgebra}
There is a Lie bracket on the complex vector space 
generated by $X_1$, $X_2$, $Y_1$, $Y_2$, $Y_3$, $Y_4$, $Z_1$, $Z_2$
such that the only non-vanishing commutators of basis elements are
$[X_1,Y_2]=Y_3$, $[X_2,Y_1]=Y_3+Z_2$, $[X_2,Y_2]=Y_4$, and $[X_2,Z_2]=Z_1$.
\end{lemma}
\begin{proof}
We have to check the Jacobi-identity:
The elements $Y_3$, $Y_4$, and $Z_1$ are in the center,
\begin{multline*}
[x_1X_1+x_2X_2+y_1Y_1+y_2Y_2+zZ_2, \\
	[x'_1X_1+x'_2X_2+y'_1Y_1+y'_2Y_2+z'Z_2,
	x''_1X_1+x''_2X_2+y''_1Y_1+y''_2Y_2+z''Z_2]] \\
= [x_1X_1+x_2X_2+y_1Y_1+y_2Y_2, (x'_2y''_1-x''_2y'_1)Z_2]
= x_2(x'_2y''_1-x''_2y'_1)Z_1,
\end{multline*}
and
\[ x'(xy''-x''y)-x''(xy'-x'y)=xx'y''-x'x''y-xx''y'+x'x''y=x(x'y''-x''y').\]
\end{proof}
The Lie algebra constructed in Lemma \ref{lem:WinkelmannLiealgebra}
is $3$-step nilpotent.
Its complex subspaces $\mf{h} = \gen{X_1,X_2}_{\C}$ and
$\mf{v} =\gen{X_1+Z_1,X_2+Z_2}_\C$ are subalgebras.
The induced $V$-action on $G/H = \exp(\gen{Y_1,Y_2,Y_3,Y_4,Z_1,Z_2}_\C) \cong \C^6$ 
is the $\C^2$-action on $\C^6$ given by $\delta$ and $\delta'$. 

Note the similarity to the construction of a pair of counterexamples to Lipsman's conjecture in
\cite{YT1}.
The smallest example of a free, affine, non-proper action of a unipotent group on some $\C^n$
is given by the $\C^2$-action on $\C^5$ generated by the two derivations 
$\delta_1 = \frac{\partial}{\partial z_1}+y_1\frac{\partial}{\partial z_2}$ and
$\delta_2 = y_1\frac{\partial}{\partial y_2}+y_2\frac{\partial}{\partial y_3}+
	y_3\frac{\partial}{\partial z_1}+\frac{\partial}{\partial z_2}$, see
\cite{YT2}, \cite{AP}.
The complex subspaces $\mf{h} =\gen{X_1,X_2}_\C$ and 
$\mf{v} = \gen{X_1+Z_1, X_2+Z_2}_\C$ are subalgebras of the $4$-step nilpotent Lie algebra
$\mf{g}$ constructed in Lemma \ref{lem:nonHausdorff}.
The induced $V$-action on $G/H \cong \exp(\gen{Y_1,Y_2,Y_3,Z_1,Z_2}_\C) \cong \C^5$ 
is the $\C^2$-action generated by $\delta_1$ and $\delta_2$, i.e., the statement of 
Corollary \ref{cor:treestep} is not true for $4$-step nilpotent Lie algebras. 
\begin{lemma}
\label{lem:nonHausdorff}
There is a Lie bracket on the vector space generated by 
$X_1$, $X_2$, $Y_1$, $Y_2$, $Y_3$, $Z_1$, $Z_2$ 
such that the only non-vanishing commutators of basis elements are
$[X_1,Y_1]=Y_2$, $[X_1,Y_2]=Y_3$, $[X_1,Y_3]=Z_2$, and $[X_2,Y_1]=Z_1$.
\end{lemma}
\begin{proof}
We verify the Jacobi-identity using the same argument as in Lemma 
\ref{lem:WinkelmannLiealgebra}.
\end{proof}
\subsection{Induced $\C$-actions on $G/H$}
Note that  $\dim H = 2 = \dim V$ in the explicit examples of the previous section.
In fact, $\dim H, \dim V \geq 2$ in all known examples for which 
$V\backslash G/H \not\cong \C^N$.
Let us look more closely at the case where one of the subgroups, say $V$, is $1$-dimensional.

Using the results of section \ref{sec:reductions} we can assume that there are 
$X_0,\ldots,X_m\in \mf{g} \setminus \mf{g}^{(l-1)}$ and $Z_0\in \mf{g}^{(l-1)}$ such
that $X_0+Z_0$ generates $\mf{v}$
and $((X_0,\ldots,X_m))$ is a Levi-Malcev basis of $\mf{h}$ if $\dim \mf{v} = 1$.
Let $\mf{g} = \mf{s} \oplus \mf{h}$ be a Levi-Malcev decomposition satisfying $Z_0\in \mf{s}$.
Recall that a $\exp\mf{s}$ is a global slice of the $H$-action on $G$
and the induced $V$-action on $G/H \cong \mf{s}$ is given by 
\[ e^{t(X_0+Z_0)}.Y = \pi_{H}(\Ad(e^{tX_0})(Y)+tZ_0) = \pi_{H}(e^{t\ad(X_0)}(Y)+tZ_0),\]
where $\pi_{H}: \mf{g} \to \mf{g}/H = \mf{s}$ 
denotes the quotient map of the $H$-action on $\mf{g}$.
Consider the linear map $\ad(X_0) : \mf{g} \to \mf{g}$.
Let $r$ be the largest $j$ such that $\ad(X_0)(\mf{g}) \subset \mf{g}^{(j)}$.
\begin{lemma}
\label{lem:dim1}
If $r\geq l-1$ or $\dim H =1$, 
then the $V\times H$-action on $G$ is proper and there is a global slice
that is algebraically isomorphic to $\C^{\dim G - \dim H - 1}$.
\end{lemma}
\begin{proof}
Since $\mf{h}$ has trivial intersection with $\mf{g}^{(l-1)}$, in both cases,
the global slice $\mf{s}$ of the $H$-action on $\mf{g}$ is $\ad(X_0)$ invariant and
contains the subspace $\gen{Z_0}_\C$.
Consequently, the induced $V$-action on $\mf{g}/H \cong \mf{s}$ is given by
\[ e^{t(X_0+Z_0)}.Y = e^{t\ad(X_0)}(Y)+tZ_0,\]
which is a free affine $\C$-action on $\mf{s}$. 

Since $\ad(X_0)(\mf{s})$ is a vector space, that does not contain $\C Z_0$ and $\C X_0$,  
there is a vector space decomposition 
$\mf{s} = \mf{s}_0\oplus  \C Z_0$ such that $\ad(X_0)(\mf{s}_0) \subset \mf{s}_0$.
Hence, $\mf{s}_0$ is a global slice of the $V$-action on $\mf{s}$.
\end{proof}

\begin{lemma}
If  $V$ is contained in a normal commutative subgroup $N$, then the $V \times H$-action
on $G$ is proper and there is a global slice isomorphic to $\C^{\dim G -\dim H -1}$.
\end{lemma}
\begin{proof}
By Lemma \ref{lem:X0maximal} we can assume that $H\subset N$, because $V/(V\cap N)$
is trivial and the $H/(H\cap N)$-action on $G/N$ admits a global slice isomorphic to
some $\C^n$. 
Let $\mf{n}=\mf{s}_1\oplus\mf{h}$ be a Levi-Malcev decomposition.
Let $S_0\subset G$ be a global slice of the $N$-action on $G$ by left-multiplication.
Then $G/H$ is isomorphic to $S_0 \times \mf{s}_1$.
The induced $V$-action on $S_0\times \mf{s}_1$ given by a derivation of the form
$\delta = \sum_j a_j(s_0) \frac{\partial}{\partial y_j}$ where the $y_j$ are 
coordinates of $\mf{s}_1$.
In other words, we obtain a family of pairs of subspaces $\Ad(s_0)(\mf{v}), \mf{h} \subset \mf{n}$
parameterized by $s_0 \in S_0$ satisfying $\Ad(s_0)(\mf{v}) \cap \mf{h} = \{ 0\}$ for all
$s_0 \in S_0$.
Because the $V$-action on $G/H \cong S_0\times S_1$ is free, there exist functions
$b_j \in \C[S_0]$ such that $\sum_j b_ja_j = 1$.
Consequently, $\delta (\sum_j b_jy_j) = \sum_j b_ja_j =1$.
\end{proof}
\begin{proposition}
\label{prop:Tiefe3}
If $\dim V =1$ and $\mf{g}^{(3)}  =\{Ê0\}$, then the $V\times H$-action on $G$ is proper and 
there is a global slice isomorphic to $\C^{\dim G -\dim H-1}$.
\end{proposition}
\begin{proof}
The Lie algebra $\mf{g}^{(1)}$ is commutative, 
since $[\mf{g}^{(j)},\mf{g}^{(k)}] \subset \mf{g}^{(j+k+1)}$.
We can assume that $\mf{h} \cap \mf{g}^{(2)} = \{ 0\}$, $\pi_0(\mf{g}) = \pi_0(\mf{v})$,
$\mf{h} = \gen{X_0,\ldots, X_m}_\C$, and $\mf{v} = \gen{X_0+Z_0}_\C$ where $Z_0$ 
is central.

The element $Z_0$ is not contained in the vector space $\mf{s}_2:= \ad(X_0)(\mf{g}^{(1)})$,
since $[X_0,Y_1] = Z_0$ would imply $\Ad(e^{-Y_1})(X_0+Z_0) = X_0$ contradicting
$\Ad(G)(\mf{v}) \cap \mf{h} =\{ 0\}$.
Let $S_2:=\exp\mf{s}_2$.
Now, the $V\times H$-action on $G/S_2$ is free if and only if the $V\times H$-action on $G$
is free, because $\ad(X_0)(\mf{g}) \cap \mf{h}$ is contained in $\mf{g}^{(1)}$ and 
$[\ad(X_0)(\mf{g}) \cap \mf{h},\mf{g}^{(1)}] = \{ 0\}$.
\end{proof}

The following examples illustrates the construction of the global slice in 
Proposition \ref{prop:Tiefe3}.
Let $G$ be the group of upper triangular $(4\times 4)$-matrices and
\[ H = \left\{ \left( \begin{smallmatrix} 1 & 0 & 0 & 0 \\ 0 & 1 & x_0 & x_1 \\
			0 & 0 & 1 & 0 \\ 0 & 0 & 0 & 1 \end{smallmatrix} \right): x_0,x_1\in \C \right\},
\,
V = \left\{ \left( \begin{smallmatrix} 1 & 0 & 0 & t \\ 0 & 1 & t & 0 \\
			0 & 0 & 1 & 0 \\ 0 & 0 & 0 & 1 \end{smallmatrix} \right): t \in \C \right\} \]
The normal commutative subgroup
\[ N = \left\{ \left( \begin{smallmatrix} 1 & 0 & y_3 & z \\ 0 & 1 & x_0 & x_1 \\
			0 & 0 & 1 & 0 \\ 0 & 0 & 0 & 1 \end{smallmatrix} \right): 
y_3,z,x_0,x_1 \in \C \right\}  \] 
contains $H$ and $V$ as subgroups.
The set
\[ S_0 = \left\{ \left( \begin{smallmatrix} 1 & y_1 & 0 & 0 \\ 0 & 1 & 0 & 0 \\
			0 & 0 & 1 & y_2 \\ 0 & 0 & 0 & 1 \end{smallmatrix} \right): y_1,y_2 \in \C \right\} \]
is a slice of the action of $N$ on $G$ by left-multiplication, i.e., 
$S_0\times N \to G$, $(s_0,n) \mapsto s_0n$ is an isomorphism.
The decomposition $\mf{n} = \{ x_0 = x_1 = 0\} \oplus \mf{h}$ defines a slice
$S_1:=\{ x_0 = x_1=0 \}$ of the $H$-action on $N$.
Using the coordinates $y_1,y_2,y_3,z$, the $V$-action on 
$G/H \cong S_0\times S_1$ is given by the derivation
$\delta = -y_1\frac{\partial}{\partial y_3}+(1-y_1y_2)\frac{\partial}{\partial z}$, because
\[Ê\left( \begin{smallmatrix} 
	1 & -y_1 & 0 & 0 \\ 0 & 1 & 0 & 0 \\ 0 & 0 & 1 & -y_2 \\ 0 & 0 & 0 & 1 \end{smallmatrix} \right)
\left( \begin{smallmatrix} 1 & 0 & 0 & t \\ 0 & 1 & t & 0 \\
			0 & 0 & 1 & 0 \\ 0 & 0 & 0 & 1 \end{smallmatrix} \right)
\left( \begin{smallmatrix} 
	1 & y_1 & 0 & 0 \\ 0 & 1 & 0 & 0 \\ 0 & 0 & 1 & y_2 \\ 0 & 0 & 0 & 1 \end{smallmatrix} \right)
= \left( \begin{smallmatrix} 1 & 0 & -y_1t & t-y_1y_2t \\ 0 & 1 & t & ty_2 
	\\ 0 & 0 & 1 & 0 \\ 0 & 0 & 0 & 1 \end{smallmatrix} \right) . \]
Hence, $\delta(z-y_2y_3) =1 $ and $S:= \{ z=y_2y_3\} \subset S_0\times S_1$ is a global
slice of the $V$-action. Note that $\log S \subset \mf{g}$ is non-linear.
Alternatively, we can construct the global slice using Lemma \ref{lem:family} with 
$G' = \{ y_2 = 0 \}$ and Lemma \ref{lem:dim1}, because $X_1 \in \mf{g}' \setminus {\mf{g}'}^{(1)}$.

Let us summarize the properties of the smallest example of a triple $H,V \subset G$ satisfying $\Ad(G)(V) \cap H = \{Êe\}$ and $\dim V = 1$ that can not be handled by the methods presented so far: 
\begin{enumerate}
\item $\mf{g}^{(3)} \ne \{ 0\}$, $\ad(X_0)(\mf{g}) \not\subset \mf{g}^{(l-1)}$ 
\item There is a basis $((X_0,\ldots,X_m))$ of $\mf{h}$ and  $Z_0\in \mf{g}^{(l-1)}$ 
	such that $X_i \in \mf{g} \setminus \mf{g}^{(l-1)}$ and $\mf{v} = \gen{X_0+Z_0}_\C$.
\item For all $Y \in \mf{g}\setminus (\pi_0^{-1}(\pi_0(\mf{h})))$,
	$\ad(Y)(\mf{h}) \not \subset \mf{h}$ and $[Y,X_0] \ne 0$.
\end{enumerate}

We finally come back to the explicit examples of the previous section that are related to
free affine $\C^2$-actions on some $\C^n$.
The quotient $\C^n/V_0$ of a subgroup $\C \cong V_0\subset \C^2$ and the induced
$\C \cong \C^2/V_0$-action on $\C^n/V$ can be easily calculated.
It would be helpful to know if these $\C$-actions, that are well understood, can
arise in our context, i.e., if they are equivalent to the $V$-action on $G/H$, where
$V$ and $H$ are subgroups of a unipotent Lie group $G$.
Unfortunately there is no method available to decide this question or 
to construct $G$, $H$ and $V$ from a given $\C$-action.

If  we choose an $\ad(X_0)$-invariant subspace $\mf{s}_0 \subset \mf{s}$ 
such that $\mf{g} = \gen{Z_0}_\C \oplus \mf{s}_0 \oplus \mf{h}$
is a Levi-Malcev decomposition as in Lemma \ref{lem:dim1}
and compatible coordinates $y_j$ of $\mf{s}_0$ and $z_0$,
then  the $\C$-action on $\mf{g}/H$ is given by a triangular derivation $\delta$ and
$\delta z_0 = 1+ P$, where $P$ is a polynomial in the coordinates $y_j$ without constant
and linear terms. 
Let the depth $d(y_j)=d$ be defined by $Y_j \in  \mf{g}^{(d-1)}\setminus \mf{g}^{(d)}$.
We define the degree $d(P):= \max_\kappa \sum_j \kappa_j d(y_j)$ of a polynomial 
$P = \sum_\kappa a_\kappa \prod_j y_j^{\kappa_j} \in \C[\mf{s}_0]$ as the usual notion
of degree weighted by the depth of the variables.
Since $[\mf{g}^{(j)},\mf{g}^{(k)}] \subset \mf{g}^{(j+k+1)}$, we have the inequality
$d(y_j) \geq d(\delta(y_j))+d(x_0)$.
We can estimate the size of $G$ assuming that a given $\C$-action 
arises as the induced $V$-action on $G/H$ using the (linear Levi-Malcev) coordinates.
For example, the quotient of the action generated by $\delta_2$ on $\C^5$ is 
$\{ z_2 = 0\} \cong \C^4$. The induced action generated by $\delta_1$ on this quotient 
is given by the derivation $-y_1^2\frac{\partial}{\partial y_2}-y_1y_2\frac{\partial}{\partial y_3}
+(1-y_1y_3)\frac{\partial}{\partial z_1}$. 
We recursively obtain $d(y_1) =1$, $d(y_2)\geq 3$, $d(y_3)\geq 5$, and $d(z_1) \geq 7$.
This means that $\mf{g}^{(6)} \ne \{ 0\}$ for the corresponding nilpotent Lie algebra $\mf{g}$.

\end{document}